\documentclass[prb,showpacs,preprintnumbers,amsmath,amssymb,twosides]{revtex4-1}

\usepackage[applemac]{inputenc}
\usepackage{amsfonts}
\usepackage{graphicx}
\usepackage{dcolumn}
\usepackage{bm}
\usepackage[english]{babel}

\newcommand{\trm}{\textrm}

\newcommand{\be}{\beqn}
\newcommand{\ee}{\eeqn}
\newcommand{\bea}{\begin{eqnarray}}
\newcommand{\eea}{\end{eqnarray}}

\newcommand{\beqn}{\begin{eqnarray}}
\newcommand{\eeqn}{\end{eqnarray}}
\newcommand{\nn}{\nonumber \\}
\newcommand{\ir}{\mathbb{IR}}
\newcommand{\irb}{\overline{\mathbb{IR}}}
\newcommand{\psia}{$\psi$-algorithm }
\newcommand{\tf}{[ f ]}

\def \A#1{\mathcal{A}_{#1}}
\usepackage{graphicx}

\begin{document}

\title{Probabilist Set Inversion using a new framework for interval arithmetic}

\author{Abdel KENOUFI}
\email{kenoufi@s-core.fr, abdel.kenoufi@uha.fr}
\affiliation{SCORE, Scientific Consulting for Research and Engineering, Strasbourg, France}
\author{Nicolas GOZE}
\email{nicolas.goze@developpement-durable.gouv.fr}
\affiliation{CETE, Centre d'Etudes Techniques de l'Equipement, Strasbourg, France}
\author{Michel GOZE}
\email{michel.goze@uha.fr}
\affiliation{LMIA, Laboratoire de Math\'ematiques, Informatique et Applications, Universit\'e de Haute-Alsace, Mulhouse, France}\date{\today}

\begin{abstract}In this paper, we present how to use an interval arithmetics framework based on free algebra construction, in order to build better defined inclusion function for interval semi-group and for its associated vector space. One introduces the $\psi$-algorithm, which performs set inversion of functions and exhibits some numerical examples developped with the \emph{python} programming langage. 
\end{abstract}

\maketitle

\section{Introduction}

Many problems dealing with real numbers can be solved with intervals computations\cite{Moore66,Moore1,Moore2,Moore3,Moore4,Sunaga}, such as set inversion, functions optimization, linear algebra. Their key stone is the construction of the inclusion function from the natural one. Some approaches have developped methods to circumvent it with using boolean inclusion tests, series or limited expansions of the natural function where the derivatives are computed at a certain point of the intervals. Nevertheless, those transfers from real functions to functions defined on intervals are not systematic and not given by a formal process. This yields to the fact that the inclusion function definition has to be adapted to each problem with the risk to miss the primitive scope.\\
We propose here a clear and simple scheme to build inclusion functions from the natural one. It is based on the arithmetic developped by one of the author in his PhD thesis\cite{nicolas}, called in this paper free-algebra based arithmetic. This article reminds first the definition and characteristics of the intervals semi-group $\ir$ and the construction of its associated vector space $\irb$. We explain how to get an associative and distributive arithmetic of intervals by embedding the vector space into a free algebra. After that, one shows from an example how to generalize the inclusion function construction for the semi-group and the vector space. This permits to present the $\psi$-algorithm in order to perform set inversion which is based on this construction and probability calculations. We ends with numerical applications examples developped in \emph{python}\cite{python} programming langage for set inversion.

\section{The Banach space $\ir$ and the associative algebra $\A{4}$}

An interval $X$ is defined as a non-empty and connected set of real numbers. It can be represented as a point in the half-plane of $\mathbb{R}^2$, $\mathcal{P}_{1}=\{(a,b)\in\mathbb{R}^2,a\le b\}$. The set $\mathcal{P}_{2}=\{(a,b)\in\mathbb{R}^2,a\ge b\}$ is the set of anti-intervals as shown on figure \ref{plan}. One writes real numbers as intervals with same bounds, $\forall a\in \mathbb{R}\ ,a\equiv [a,a]$. 
We denote by $\ir=\mathcal{P}_{1}$ the set of intervals of $\mathbb{R}$. 
The semantical arithmetic operations on intervals, called \textbf{\emph{Minkwoski or classical operations}},  are defined such that the
result of the corresponding operation on elements belonging to operand
intervals belongs to the resulting interval. That is, if $\diamond $ denotes
one of the usual operations $+,-,\ast,/ $, we have,
if $X$ and $\ Y$ are bounded intervals of $\mathbb{R}$,
\begin{equation*}
X\diamond Y=\{x\diamond y\text{ }/\text{ }x\in X,\text{ }y\in Y\},
\end{equation*}
But $\ir$ is only a semi-group for the addition and no substraction can be defined for it. 

\subsection{Modal interval arithmetic}

An algebraic extension of the classical interval arithmetic, called generalized interval arithmetic\cite{Sunaga,Moore66} has been proposed first by M. Warmus\cite{Warmus1,Warmus2}. It has been followed in the seventies by H. -J. Ortolf and E. Kaucher \cite{Ortolf,Kaucher1,Kaucher2,Kaucher3,Kaucher4}. In this former interval arithmetic, the intervals form a group with respect to addition and a complete lattice with respect to inclusion. In order to adapt it to semantic problems, Gardenes et al. have developed an approach called modal interval arithmetics\cite{Gardenes1,Gardenes2,Gardenes3,Gardenes4,Gardenes5,Gardenes6}.
S. Markov and others investigate the relation between generalized intervals operations and Minkowski operations for classic intervals and propose the so-called directed interval arithmetic, in which Kaucher's generalized  intervals can be viewed as classic intervals plus direction, hence the name directed interval arithmetic\cite{Markov1,Markov2}. In this arithmetic framework, proper and improper intervals are considered as intervals with sign\cite{Markov3}. Interesting relations and developements for proper and improper intervals arithmetic and for applications can be found in litterature\cite{Dimitrova,Popova1,Popova2}.\\
Our approach\cite{nicolas,rc1}, that we remind below in this artice, is similar to the previous ones in the sense that intervals are extended to generalized intervals. However we use a construction which is more canonical and based on the semi-group completion into a a group, which permits then to build the associated real vector space, to get an analogy with directed interval arithmetic and to perform differential calculus and linear algebra naturally.

\subsection{New interval arithmetic}

Let $\irb$ be the $\ir$ completed group. An element of $\irb$ is written $\overline{(A,B)}$ where $A,\ B\in \mathbb{IR}$ and the equivalence relation being $(A,B)\sim(A',B')$ if and only if $A+B'=B+A'$. For example $([3,4],[1,4])\sim([2,0],[0,0])\equiv([2,0],0)$. Any equivalence class can be written as $\overline{(X,0)}$ or $\overline{(0,X)}$ or $\overline{(a,a)}$ with $X\in \ir$ or $a\in \mathbb{R}$. We denote $\mathcal{X}\setminus \mathcal{Y}, \mathcal{X},\  \mathcal{Y}\in \irb$, the difference into the abelian group $\irb$. The external multiplication
\beqn
\mathbb{R}\times \irb \longrightarrow \irb
\eeqn
given by 
\beqn
a\cdot \overline{(X,0)}&=&\overline{(a\cdot X,0)}\  \trm{if}\  a\ge 0\\
&=&\overline{(0,-a\cdot X)}\  \trm{if}\  a\le 0\\
\eeqn
and 
\beqn
a\cdot \overline{(0,X)}&=&\overline{(0,a\cdot X)}\  \trm{if}\  a\ge 0\\
&=&\overline{(-a\cdot X,0)}\  \trm{if}\  a\le 0\\
\eeqn
provides $\irb$ with a 2-dimensional $\mathbb{R}$ vector space structure. A basis is given by the family $\{\overline{([1,1],0)},\overline{([0,1],0)}\}$.
The following map
\beqn
||.|| : \irb &\longrightarrow& \mathbb{R}^+\\
\overline{(X,0)}&\mapsto& l(X)+|c(X)|\\
\overline{(0,X)}&\mapsto&l(X)+|c(X)|\\
\eeqn
with respectively $l(X)$ and $c(X)$ the width and the center of the interval $X$, is a norm. It is easy to prove\cite{nicolas,rc1} that $\irb$ is isomorphic to $\mathbb{R}^2$ which is complete, this yields to the  fact that this norm endows $\irb$ with a Banach space structure.
Moreover, one can define an order relation $<$ on $\ir$ which is compatible with this norm : $X,Y\in \mathbb{IR}, X<Y$ if
\begin{equation}
 X\not\subset Y \trm{ and } c(X)<c(X) 
\end{equation}
or
\begin{equation}
 X\subset Y \trm{ and } l(X)<l(Y)
\end{equation}
Due to the fact that every intervals can be compared to each others, it yields that this partial order relation can be extended to a total one.\\
In order to define a natural multiplication on $\irb$, we consider the following 4-dimensional associative algebra $\A{4}$ with basis :
\begin{equation*} 
\left\{
\begin{array}{l}
e_{1}=(1,1), \\
e_{2}=(0,1), \\
e_{3}=(-1,0), \\
e_{4}=(-1,-1).%
\end{array}%
\right.
\end{equation*}%
and product :
\begin{equation}
\label{minkowski}
\begin{tabular}{|l|l|l|l|l|}
\hline
& $e_{1}$ & $e_{2}$ & $e_{3}$ & $e_{4}$ \\ \hline
$e_{1}$ & $e_{1}$ & $e_{2}$ & $e_{3}$ & $e_{4}$ \\ \hline
$e_{2}$ & $e_{2}$ & $e_{2}$ & $e_{3}$ & $e_{3}$ \\ \hline
$e_{3}$ & $e_{3}$ & $e_{3}$ & $e_{2}$ & $e_{2}$ \\ \hline
$e_{4}$ & $e_{4}$ & $e_{3}$ & $e_{2}$ & $e_{1}$ \\ \hline
\end{tabular}%
\end{equation}%
\begin{figure}[!h]
\label{plan}
\centering
\includegraphics[width=20cm,height=10cm]{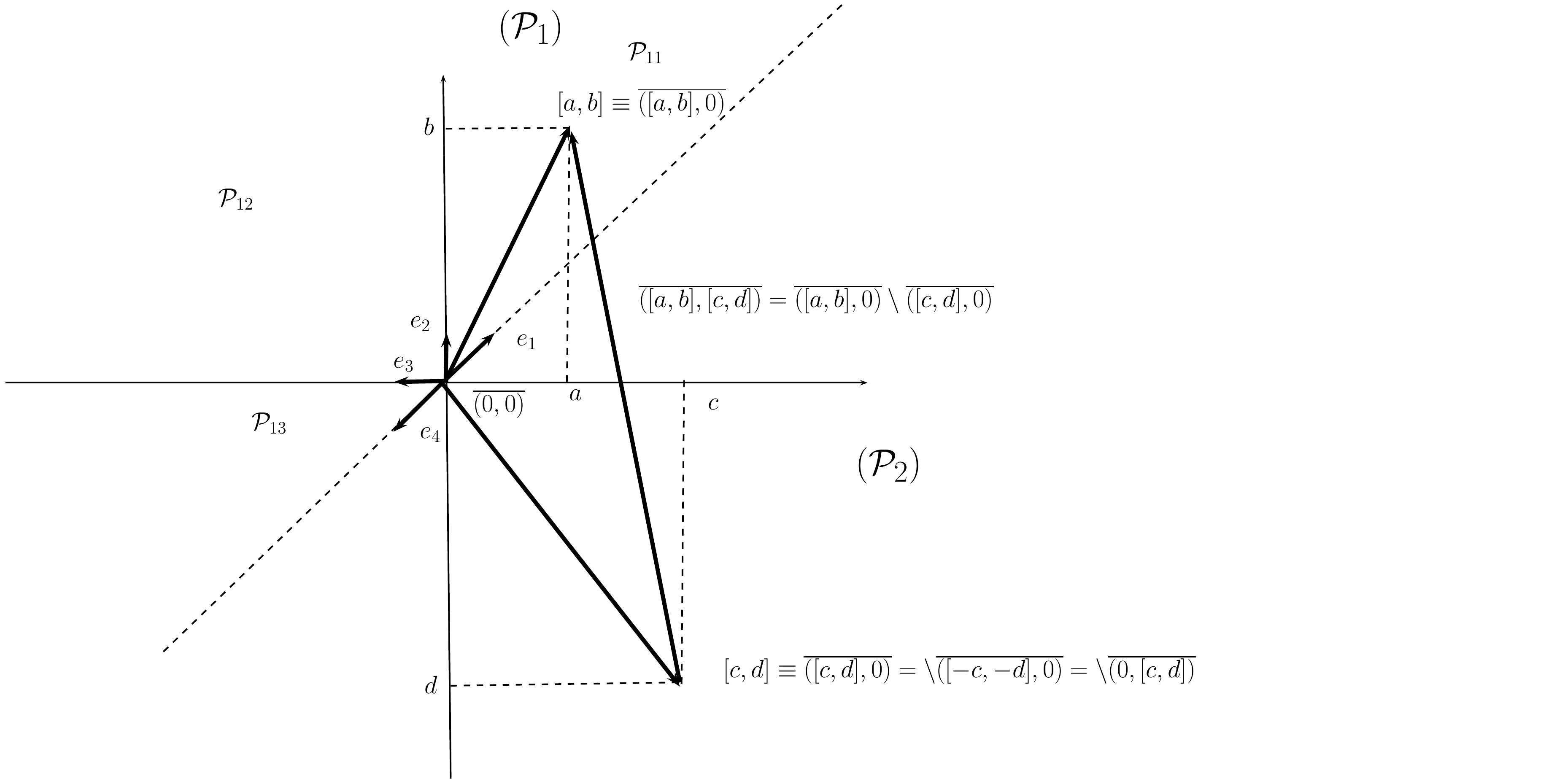}
\caption{Intervals and anti-intervals can be respectivelly represented as points in half planes $\mathcal{P}_{1}$ and $\mathcal{P}_{2}$. This implies naturally the representation of interval vectors from $\irb$. }
\label{plan}
\end{figure}

Those vectors represent the intervals $[1,1],[0,1],[-1,0]$ and $[-1,-1]$ (figure \ref{plan}). The product table is obtained simply with the mean of Minkoswki product. For example $e_2\cdot e_3=[0,1]\cdot[-1,0]=[\min(0,-1),\max(0,-1)]=[-1,0]=e_3$.
One can define the sectors $\mathcal{P}_{1j}$, with $j=1,2,3$ such as 
\beqn
{j=1 :\ }(a,b)&=&ae_1+(b-a)e_2,\ \trm{if }a\ge 0 .\\
{j=2 :\ }(a,b)&=&-ae_3+be_2, \trm{if }a\le 0\ \trm{and }b\ge 0.\\
{j=3 :\ }(a,b)&=&-be_4+(b-a)e_3, \trm{if } b<0.
\eeqn

\begin{figure}[!h]
\centering
\includegraphics[width=12cm,height=5cm]{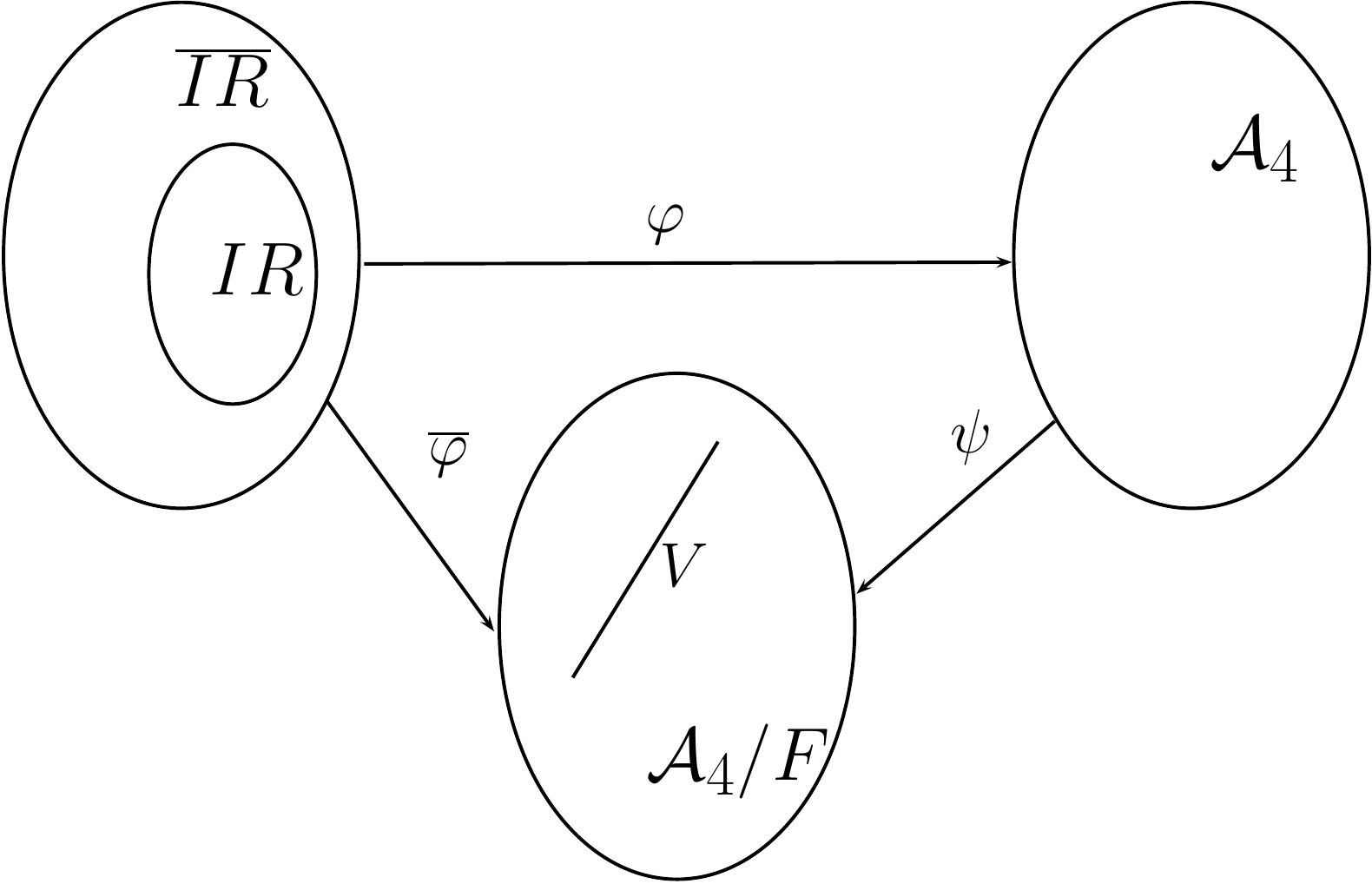}
\caption{Interval embedding in the algebra $\mathcal{A}_4$.}
\label{a4}
\end{figure}

Thus any interval $X$ (i.e any $\overline{(X,0)}\in\irb$) is represented by positive components vector of $\A{4}$.
Let's denote this embedding $\varphi : \ir \longrightarrow \A{4}$. The representation of $\overline{(0,X)}$ is then given by $\varphi(\overline{(0,X)})=-\varphi(\overline{(X,0)})$.
Let $F$ be the linear subspace of $\A{4}$ spanned by $e_1-e_2+e_3$ and $e_1+e_4$. Any vector $\sum_i \alpha_i e_i\in \A{4}$ is equivalent in $\A{4}/F$ to a vector with positive components if $\alpha_2+\alpha_3\ge 0$. In this case, such a vector corresponds to an interval of $\ir$ and one has the bijection $\overline\varphi : \irb \longrightarrow V$ where $V$ is the convex part of $\A{4}/F$ corresponding to $\alpha_2+\alpha_3\ge 0$. This defines a distributive product in $\ir$ by
\bea
X\cdot Y={\overline{\varphi}}^{-1}(\varphi'(X)\cdot\varphi'(Y))
\eea
where $\varphi '=\psi\circ\varphi$ and $\psi$ is the canonical embedding from $\A{4}$ to  $\A{4}/F$ as shown on Fig.(\ref{a4}).
Division between intervals can also be defined with solving $X\cdot Y=(1,0,0,0)$ in $\A{4}$ or in an isomorphic algebra.

One can remark that :
\begin{itemize}
\item Distributivity is satisfied for multiplication and division according addition and substraction
\item Usually, one has $X\diamond Y \subseteq X\cdot Y$ where $X\diamond Y$ is the classical Minkowski product. 
More precisely, within this algebra, if $X$ and $Y$ are not in the same sector $\mathcal{P}_{1j}$, or are in the same sector $\mathcal{P}_{11}$ or  $\mathcal{P}_{13}$, one gets the Minkowski product. If they are belonging both to $\mathcal{P}_{12}$, the result is greater than the Minkowski one. If one needs a result $X\cdot Y$ closer to $X\diamond Y$ or even equal to, one can embed $\ir$ in higher dimension algebras, such as $\A{5}, \A{6},\A{7}$, and so on. This former one is sufficient to recover the distributivity and equality to Minkowski product\cite{nicolas}. 
\item It can be naturally extended to $\irb$.
\item The product satisfies the monotony property.
\item With the canonical injection $X\mapsto \overline{(X,0)}$ from $\ir$ to $\irb$, one can consider that $\ir\subset\irb$.

\end{itemize}
\section{Inclusion functions}
It is necessary for some problems to extend the definition of a function defined for real numbers $f : \mathbb{R}\rightarrow \mathbb{R}$ to function defined for intervals $\tf : \ir \rightarrow \ir$ such as  $\tf([a,a])=f(a)$ for any $a\in\mathbb{R}$. It will be convenient to have the same formal expression for $f$ and $\tf$. Usually the lack of distributivity in Minkowski arithmetics doesn't give the possibility to get the same formal expressions. But with the arithmetic we have presented, there is no data dependency anymore and one can define easily inclusion functions from the natural one. For example, let's extend to intervals the real functions $f_1(x)=x^2-2x+1$, $f_2(x)=(x-1)^2$, $f_3(x)=x(x-2)+1$. 
Usually, with the Minkwoski operations, the three expressions of this same function for the interval $X=[3,4]$ are $\tf_1(X)=[2,11]$, $\tf_2(X)=[4,9]$ and $\tf_3(X)=[6,12]$. Data dependancy occurs when the variable appears more than once in the function expression. The deep reason of that is the lack of distributivity in Minkowski arithmetics. But within the arithmetic developped in $\A{4}$ or higher dimension free algebras\cite{nicolas,rc1}, this problem vanishes. For example : with $X=[3,4]$ and since $X \in \mathcal{P}_{11}$, 
\be 
\varphi(X)=(3,4-3,0,0)=(3,1,0,0)=3e_1+e_2.
\ee
Since $e_1\equiv 1$, and with means of product table (\ref{minkowski}) one has
\be 
\varphi(\tf_1(X))&=&(3e_1+e_2)^2-2(3e_1+e_2)+1\nn
&=&9e_1^2+2\cdot 3e_1e_2+e_2^2-2\cdot 3e_1-2e_2+1\nn
&=&9e_1+6e_2+e_2-6e_1-2e_2+e_1=4e_1+5e_2\nn
&=&\varphi([4,9]),
\ee

\be 
\varphi(\tf_2(X))&=&(3e_1+e_2-1)^2\nn
&=&(2e_1+e_2)^2\nn
&=&4e_1^2+4e_1e_2+e_2^2\nn
&=&4e_1+4e_2+e_2\nn
&=&4e_1+5e_2\nn
&=&\varphi([4,9]),
\ee
and

\be 
\varphi(\tf_3(X))&=&(3e_1+e_2)\cdot(3e_1+e_2-2)+1\nn
&=&9e_1+3e_1e_2-6e_1+3e1e_2+e_2^2-2e_2+e_1\nn
&=&4e_1+3e_2+3e_2+e_2-2e_2\nn
&=&4e_1+5e_2\nn
&=&\varphi([4,9]).
\ee

Thus, $\tf_1(X)=\tf_2(X)=\tf_3(X)=[4,9]$.\\
On the other hand, the construction of the inclusion function depends on the type of problem one deals with. If one aims to perform set inversion for example, it has to be done in the semi-group $\ir$ due to semantic definition of substraction. But, the substraction is not defined in $\ir$. This problem can be circumvented by replacing it with an addition and a multiplication with the interval $e_4=[-1,-1]$. This maintains the associativity and distributivity of arithmetics and permits to introduce a semantical substraction. For example : if $f(x)=x^2-x=x(x-1)$ for real numbers, one defines $\tf(X)=X^2+e_4\cdot X$. One reminds the product $[-1,-1]\cdot [a,b]$ is equal to $[-b,-a]$. Due to the fact that the arithmetics is now associative and distributive, one doesn't have data dependancy anymore and  $\tf(X)=X^2+e_4\cdot X=X\cdot (X+e_4)$. The last term corresponds to the transfer of $x(x-1)$. Division can be transfered to the semi-group in the same way by replacing $\frac{1}{x}=x^{-1}$ with $X^{e_4}$.    \\
Taylor polynomial expansions, differential calculus and linear algebra operations are defined only in a vector space. Therefore the transfer for the vector space is done directly. This permits to get infinitesimal intervals with the substraction and to compute derivatives. This is of course not allowed and not possible into the semi-group. From $\ir$ to the vector space $\irb$, $f:x\mapsto -x$ is transfered to $\tf:X\equiv \overline{(X,0)}\mapsto \backslash \overline{(X,0)}\equiv \backslash X$. This means that  $[a,b]$ substraction is the anti-interval $[-a,-b]$ addition. One of the most important consequence is that it is possible to transfer some functions directly to the intervals. For example, it is easy to prove analytically in $\irb$ that $[\exp](\overline{([a,b],0)})=\overline{([\exp(a),\exp(b)],0)}$ with means of Taylor expansion.

\section{Probabilist Set Inversion : \psia}
\subsection{Flowchart}
\begin{figure}[!h]
\centering
\includegraphics[width=10cm,height=5cm]{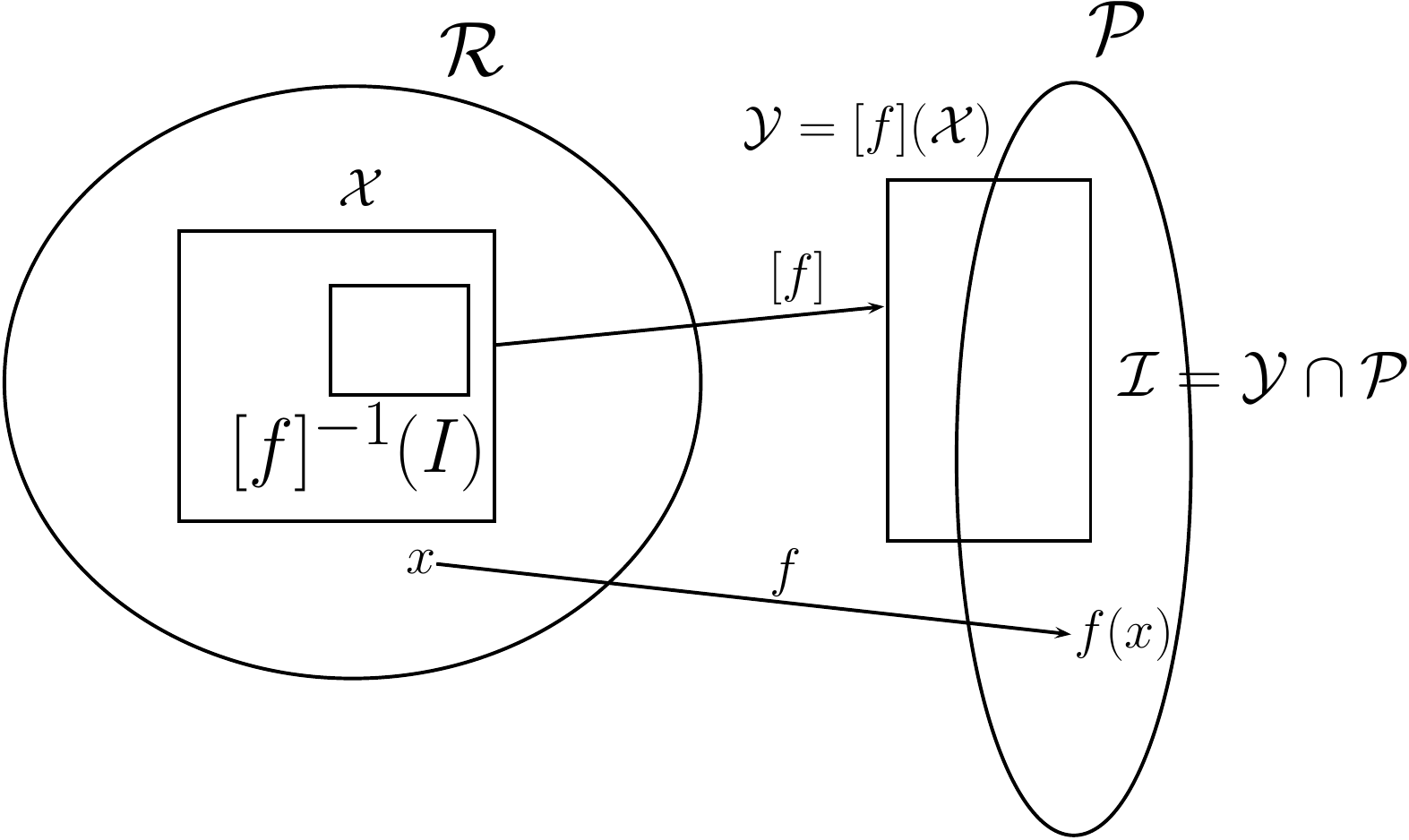}
\caption{Set inversion modelization with intervals}
\label{psi}
\end{figure}

One presents an efficient set inversion method whose flowchart is very simple. One of the powerful application of interval calculus is the set inversion of a real-valued function defined on real numbers. It finds a huge number of applications in automation, optimization, resolution of algebraic and differential equations. The mathematical modelization of this problem is the following as shown on Fig.(\ref{psi}) : let's note $f:\mathbb{R}^n\mapsto\mathbb{R}^p$ a function for a physical system, which is required to be surjective only, $\mathcal{R}\subset\mathbb{R}^n$ the set of adjustements, and $\mathcal{P}\subset\mathbb{R}^p$ the set of performance of a system. Set inversion consists of the computation of $\mathcal{S}=f^{-1}(\mathcal{P})\cap\mathcal{R}$. For semantical reasons, one has to perform this set inversion within the semi-group $\ir$.

\begin{figure}[!h]
\centering
\includegraphics[width=12cm,height=8cm]{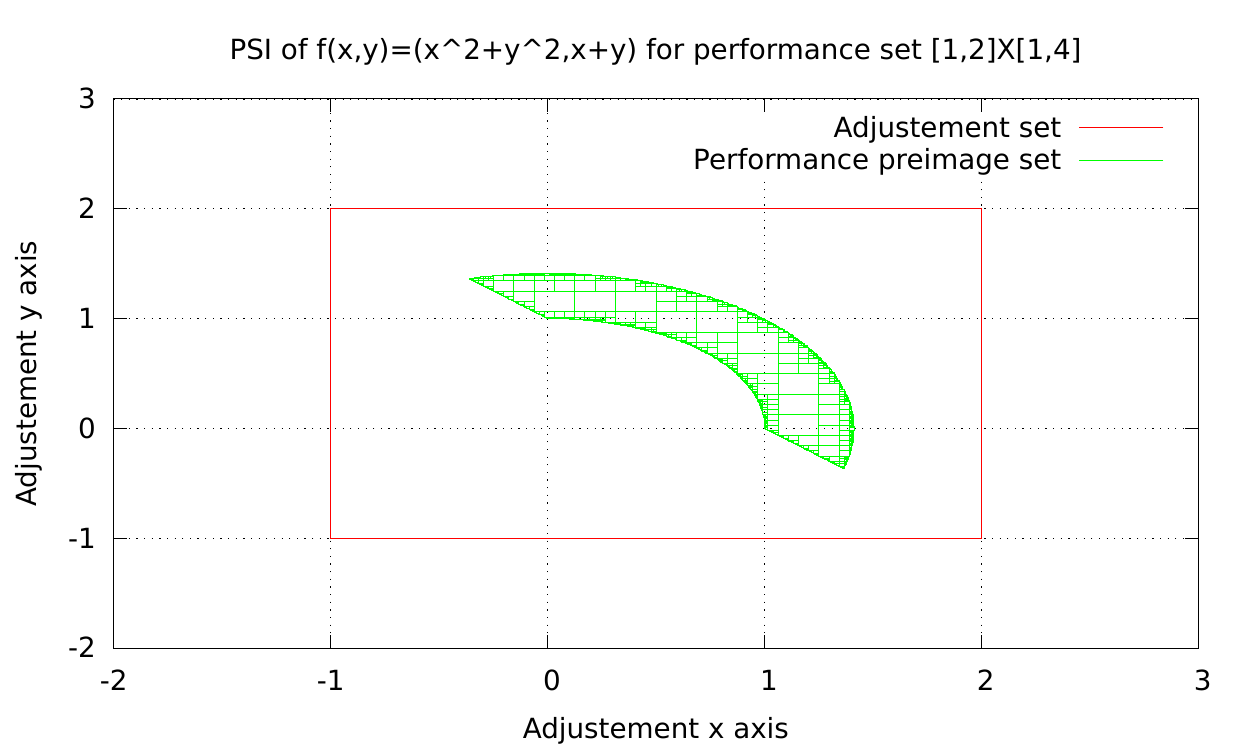}
\caption{\psia for $f_1$}
\label{psi1}
\end{figure}

\begin{figure}[!h]
\centering
\includegraphics[width=12cm,height=8cm]{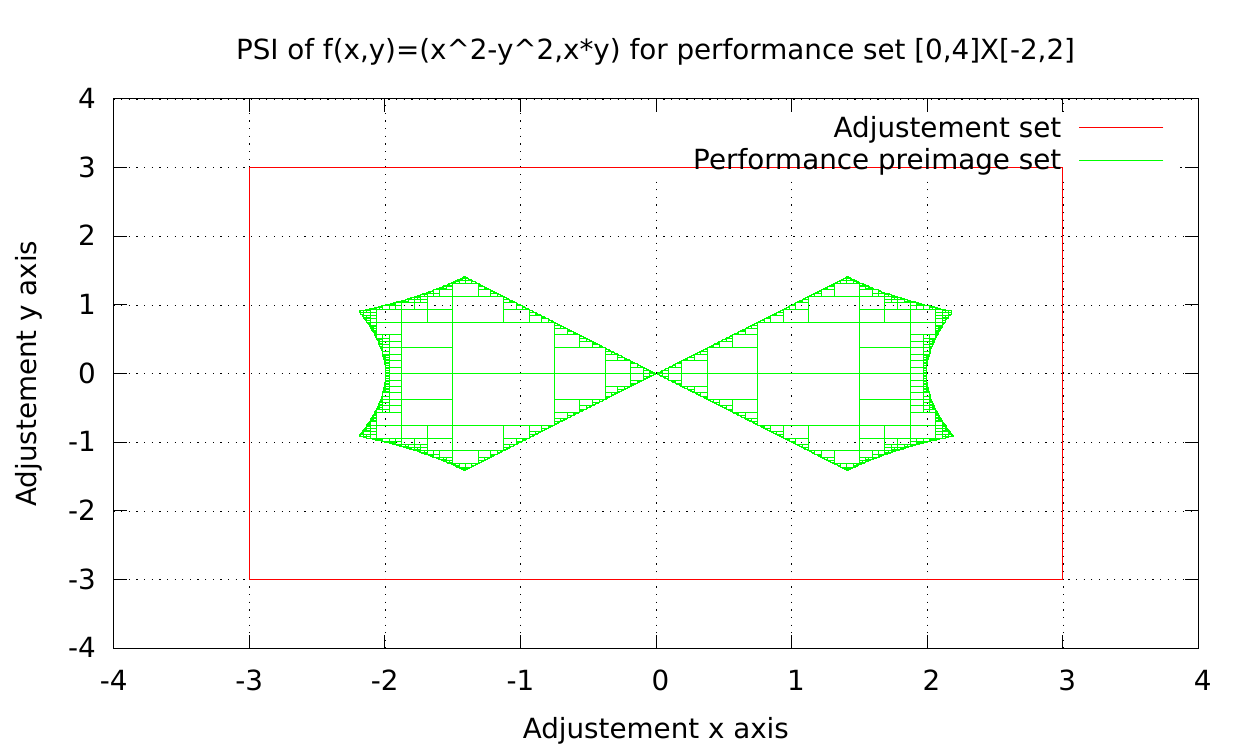}
\caption{\psia for $f_2$}
\label{psi2}
\end{figure}

\begin{figure}[!h]
\centering
\includegraphics[width=12cm,height=8cm]{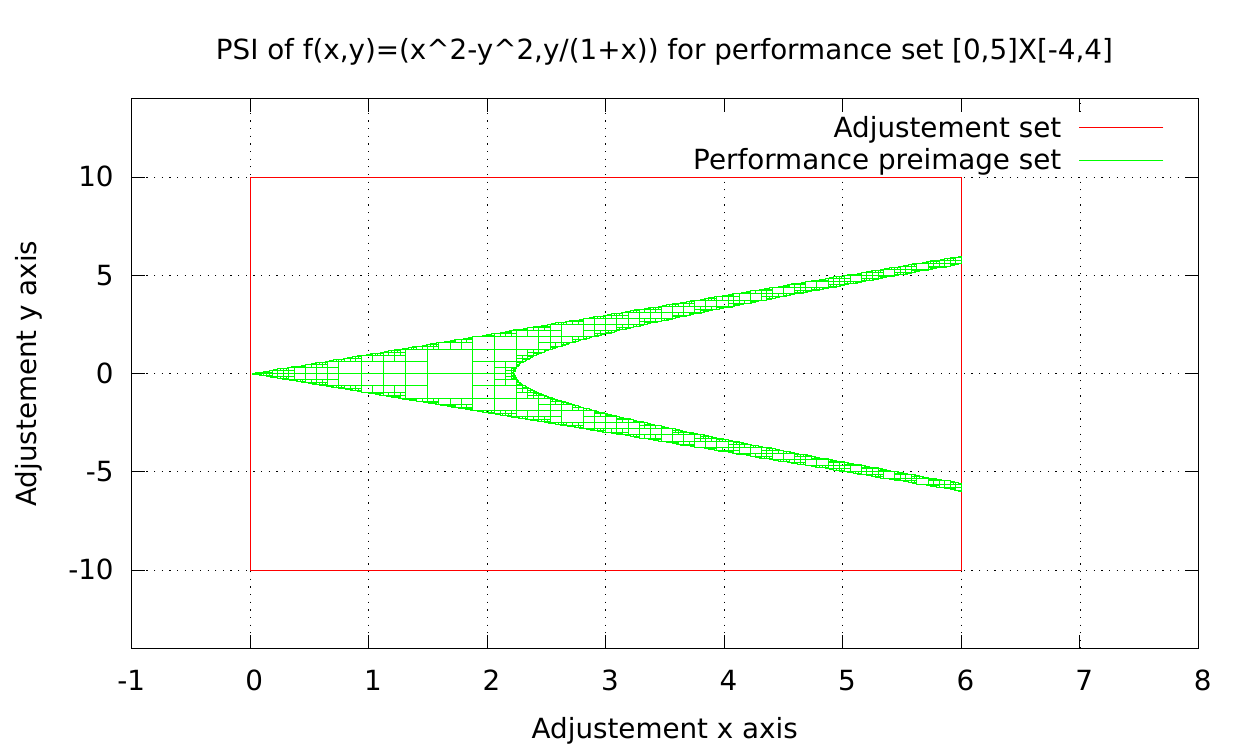}
\caption{\psia for $f_3$}
\label{psi3}
\end{figure}

Some interesting and powerful methods using intervals have been developped those last years, especially SIVIA\cite{Jaulin}, Set Inversion \emph{Via} Interval Analysis. But the inclusion function being not well defined in the semi-group with the Minkowski arithmetic, SIVIA uses boolean inclusion tests and finds accepted, rejected and "uncertain" domains.\\
With the algebraic arithmetic, one doesn't need boolean tests since the inclusion functions are well-defined. Thus, we propose the \psia (Probabilist Set Inversion) inspired from SIVIA but without boolean tests and with a conditional probability calculation and domain bissections. This yields to accepted or rejected domains only.
We are interested to compute the following conditional probability 
\beqn
p(\mathcal{X})&=&p([f](\mathcal{X})\subset\mathcal{P}\ | \  f(x)\in[f](\mathcal{X}),\ \forall x\in \mathcal{X})\nn
&=&\frac{mes([f](\mathcal{X})\cap\mathcal{P})}{mes([f](\mathcal{X}))}\nn
&=&\frac{mes(\mathcal{Y}\cap\mathcal{P})}{mes(\mathcal{Y})}=\frac{mes(\mathcal{I})}{mes(\mathcal{Y})}\nn
\eeqn
where $mes$ is the Lebesgue measure in $\mathbb{R}^p$ (length, surface, ...).
If this probability equals $1$ then the set is added to the list of solutions. If it is zero the set is rejected and removed from the list of interval candidates. If the probability is such as $p(\mathcal{X})\in ]0,1[$, then $\mathcal{X}$ is bissected and \psia applies the same procedure recursively for the resulting intervals until the size is lower than a fixed size resolution of the intervals or until the sets are accepted or rejected. Since \psia creates sequences of decreasing intervals which are compact sets, it is obvious that \psia converges to fixed points probabilities which are simply $0$ and $1$. In fact, we consider a sequence of compact sets $\{K_n\}_{n\in \mathbb{N}}$ satisfying $d(K_n)< d(K_{n+1})$ where $d$ is the diameter of the compact set. If $K_n\cap K_{m} = \emptyset$, for any $n\ne m$, the sequence is convergent and the limit is the empty set. Then, one has just to consider the sequence of bissected sets.

\begin{figure}[!h]
\centering
\includegraphics[width=12cm,height=8cm]{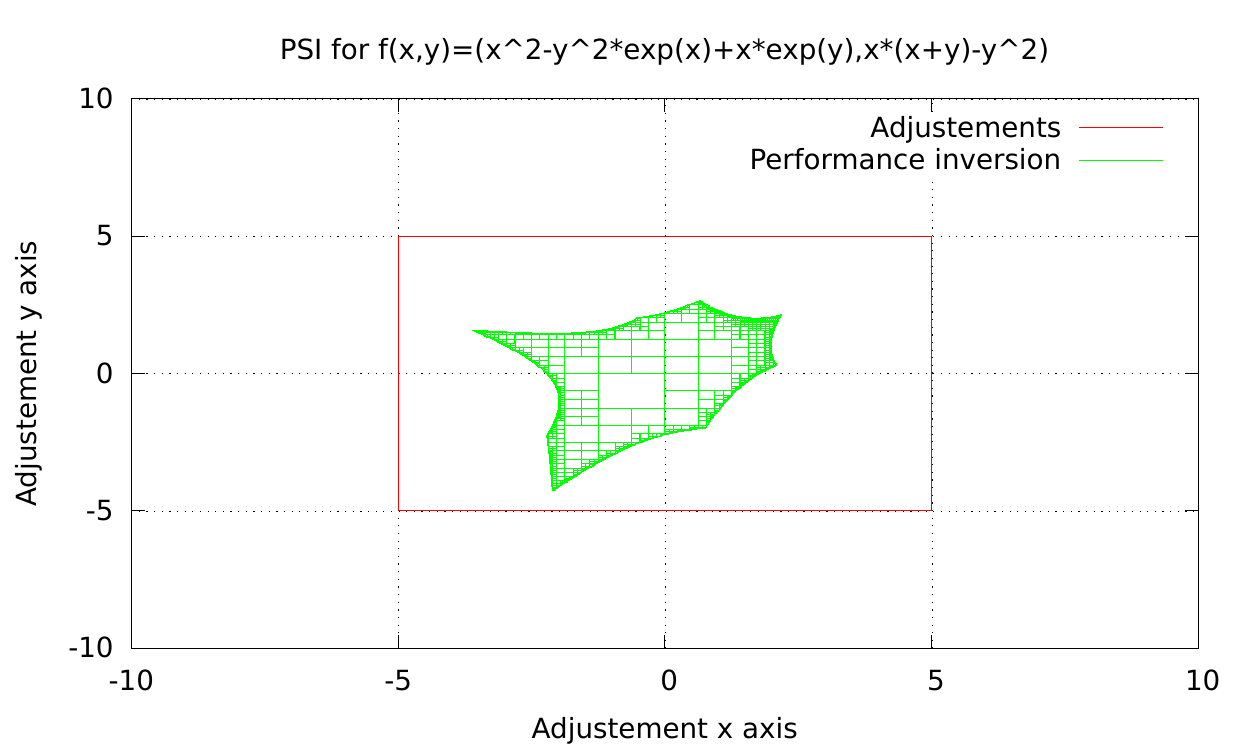}
\caption{\psia for $f_4$}
\label{psi4}
\end{figure}

\begin{figure}[!h]
\centering
\includegraphics[width=12cm,height=8cm]{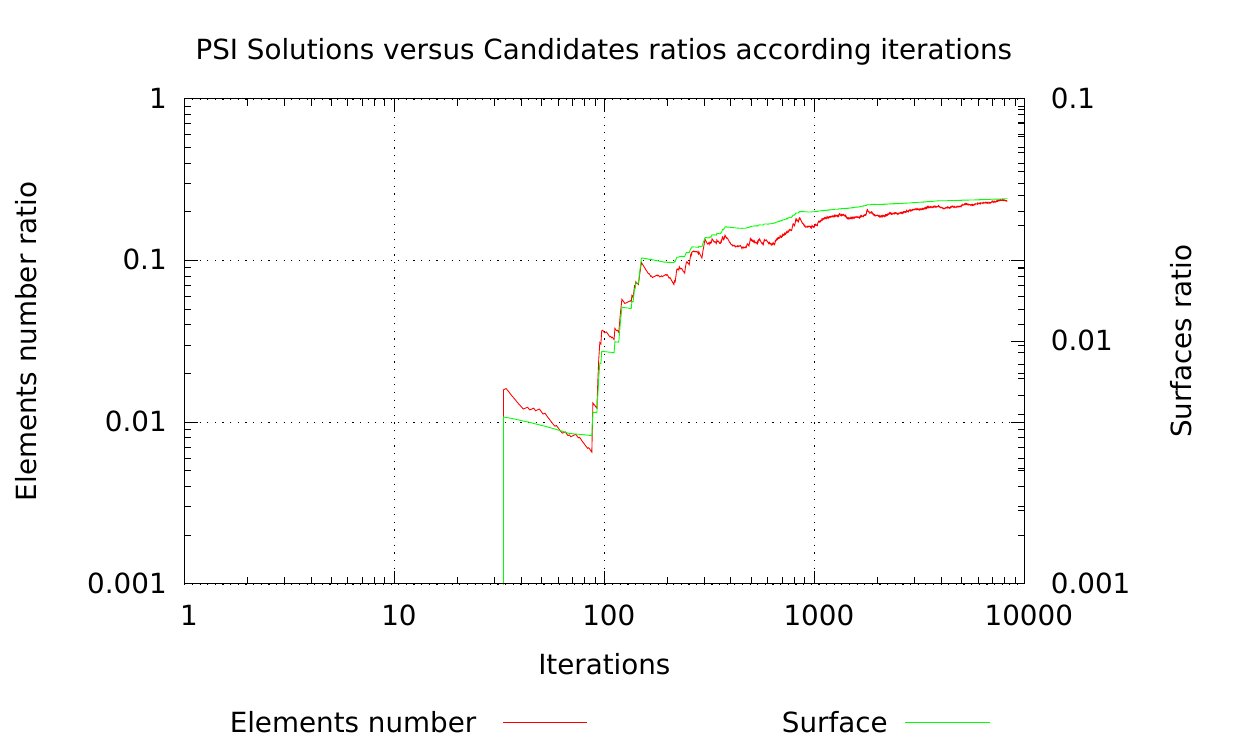}
\caption{Accepted and candidate intervals total ratios for numbers and surfaces according iterations in \psia for $f_4$ with $10^{-4}$ spatial resolution.}
\label{iterations}
\end{figure}

\subsection{Numerical applications}

We have developped a numerical library for \emph{python} environnement\cite{python} called yet \emph{interval\_lib}. It is a pure numerical implementation performing the basic arithmetics presented above \cite{nicolas,rc1}. This library aimes to give simple and optimized routines to perform interval calculations based on the algebraic arithmetics. One gives in this section some numerical application examples of \psia in order to illustrate how it can treat usual inversion problems and build well-defined inclusion functions.\\
Let's define the non-linear functions $f_i\ :\ \mathbb{R}^2 \rightarrow \mathbb{R}^2$ with respectively adjustements and performances sets $\mathcal{R}_i$, $ \mathcal{P}_i$ :
\be
f_1(x,y)&=&(x^2+y^2,x+y),\ \mathcal{R}_1=[-1,2]^2,\ \mathcal{P}_1=[1,2]\times [1,4]\nn
f_2(x,y)&=&(x^2-y^2,x\cdot y),\ \mathcal{R}_2=[-3,3]^2,\ \mathcal{P}_2=[0,4]\times [-2,2]\nn
f_3(x,y)&=&(x^2-y^2,\frac{y}{1+x}),\ \mathcal{R}_3=[0,6]\times [0,10],\ \mathcal{P}_3=[0,5]\times [-4,4].\nn
\ee
Those examples have been chosen to give examples of addition, substraction, product and division transfers from $\mathbb{R}$ to $\ir$, and to exhibit the difference between the usual Minkowski arithmetic and the algebraic  one\cite{Jaulin2,nicolas,rc1}. The calculations with the \psia are shown on figures \ref{psi1}, \ref{psi2}, \ref{psi3} and \ref{psi4}. The convergence to $0$ or $1$ probabilities only, shows that inclusion functions are well constructed and that the interval arithmetic is robust.
The last example 
\be
f_4(x,y)=(x^2-y^2\cdot\exp(x)+x\cdot\exp(y),x\cdot(x+y)-y^2),\ \mathcal{R}_4=\mathcal{P}_4=[-5,5]^2\nn
\ee 
presented on figure \ref{psi4} shows clearly that the \psia implemented in the algebraic arithmetic we use is not data dependant. The variables appear more than once in the formal expression of the function $f_4$. The CPU time for this inversion is about $255$ seconds on a simple $1.67$ Ghz Intel processor for a spatial surface resolution of $10^{-4}$.\\
In figure \ref{iterations}, one computes for this resolution the number of solution intervals according to candidates ones and the  total surfaces ratio of the same intervals. The strong correlation between the two curves caracterizes the fact that the measures computations are well-defined because the probabilities are very close to the frequencies computed with \psia independantly of the initial adjusments set. This confirms that the free algebra based arithmetic framework is relevant even for a semi-group $\ir$ calculations example. Due to the bissection, the algorithm computational complexity is exponential according to the iterations $N$, and it is not improved compared to SIVIA one. In our scheme, computational time is defined as 
\be
T_{comp}=\mathcal{O}(N)=k\cdot 2^N.
\ee 
However, if the native function is differentiable on $\mathbb{R}^n$, it is possible to define an adaptative mesh, with bissection spanned only in the space directions where the derivative magnitude is larger than a certain fixed value, because it is not useful to bissect in flat directions. This will obviously decrease the computational complexity constant $k$.

\section{Conclusion}
We have presented in this article an heuristic way to transfer real functions to inclusion ones depending on the space needed, semi-group or vector space of intervals. This permits to define a simple but very efficient algorithm for set inversion inspired from SIVIA\cite{Jaulin}, the \psia, which uses a free algebra based arithmetic \cite{nicolas,rc1} and probability calculations. The convergence of this algorithm is assured for linear and non-linear problems. It offers several possibilities of applications, such as solving algebraic equations, differential equations, probability law of random variables calculations (discrete or continuous), topological analysis, numerical Lebesgue integrals computations, data analysis such as principal components analysis, and parameters identification.

\begin{acknowledgements}
We thank Michel Gondran, Luc Jaulin and Jean-Fran\c cois Osselin for useful and interesting discussions.
\end{acknowledgements}


\end{document}